# Linear Parameter Varying Power Regulation of Variable Speed Pitch Manipulated Wind Turbine in the Full Load Regime

T. SHAQARIN[1], MAHMOUD M. S. AL-SUOD[2]
[1]Department of Mechanical Engineering, TafilaTechnical University, Tafila 11660, JORDAN
[2]Department of Electrical Power Engineering and Mechatronics, TafilaTechnical University, Tafila 11660, JORDAN

*Abstract:* - In a wind energy conversion system (WECS), changing the pitch angle of the wind turbine blades is a typical practice to regulate the electrical power generation in the full-load regime. Due to the turbulent nature of the wind and the large variations of the mean wind speed during the day, the rotary elements of the WECS are subjected to significant mechanical stresses and fatigue, resulting in conceivably mechanical failures and higher maintenance costs. Consequently, it is imperative to design a control system capable of handling continuous wind changes. In this work, Linear Parameter Varying (LPV) $H_\infty$ controller is used to cope with wind variations and turbulent winds with a turbulence intensity greater than ± 10%. The proposed controller is designed to regulate the rotational rotor speed and generator torque, thus, regulating the output power via pitch angle manipulations. In addition, a PI-Fuzzy control system is designed to be compared with the proposed control system. The closed-loop simulations of both controllers established the robustness and stability of the suggested LPV controller under large wind velocity variations, with minute power fluctuations compared to the PI-Fuzzy controller. The results show that in the presence of turbulent wind speed variations, the proposed LPV controller achieves improved transient and steady-state performance along with reduced mechanical loads in the above-rated wind speed region.

## 1 Introduction

During the past decades, the usage of wind turbine plants increased and became more competitive among other renewable energy forms. The current trend is toward the advancement of green energy production. On the contrary, more restrictions are enforced to reduce the energy produced from conventional generation sources, which produce greenhouse gases that lead to global warming. Recently, the wide adoption of wind farms in the United States positively affects greenhouse emissions and water consumption. For instance, in 2018, the United States' wind capacity of 96 GW annually lessened $CO^2$ emissions and water consumption by about 200 million metric tons, and 95 billion gallons, respectively, [1].

Due to the turbulent nature of the wind and the large variations of the mean wind speed during the day, the rotary elements of the WECS are subjected to significant mechanical stresses and fatigue, resulting in conceivably mechanical failures and higher maintenance costs.

As a result, mechanical failures are unavoidable in dynamic and turbulent wind speed conditions unless proper power regation control systems are used, such as pitching controllers. Despite decreasing mechanical stresses on the wind turbine's rotary elements, a sophisticated control system helps in stabilizing, maximizing, and restricting the generated power at above-rated wind speeds, [2].

In wind energy generation systems, the effect of turbulence variations and noise created by turbines may cause fluctuations and a reduction in power output. Wind turbines are designed and analyzed by modeling and simulating the turbines using various software and hardware techniques. The wind turbine model in WECS was developed by Manyonge *et al.*, [3], via examining the power coefficient parameter needed to understand the wind turbine dynamics over its operational regime, which contributes to controlling the performance of wind turbines.

The work presented by Taher *et al.*, [2], introduced a Linear Quadratic Gaussian (LQG) gain-scheduling controller to cope with the variable wind velocity in the WECS. They introduced gain-scheduling controllers (GSC), which aimed to manipulate the blades' pitch angle, regulate the electrical torque, and maintain the rotor speed constant at their nominal values at the full load region. Petru and Thiringer, [4], presented in their work a dynamic model for a wind turbine in both fixed-speed and stall-regulated systems. The main objective of his work is to evaluate the dynamic power quality effect on the electrical grid.

Sabzevari *et al.*, [5], proposed a maximum power point tracking (MPPT) approach emanating from a neural network that is trained offline using particle swarm optimization (PSO). The proposed MPPT estimated the wind speed to adapt the fuzzy-PI controller. The adaptive controller manipulated the boost converter duty cycle for the driven permanent magnet synchronous generator (PMSG). Macêdo and Mota, [6], presented a comprehensive description of the wind turbine system equipped with an asynchronous induction generator. They implemented a controller using a Fuzzy control system by manipulating the pitch angle. The main objective of their work was to reduce the fluctuations in the generator output power. Salmi *et al.*, [7], designed an optimal backstepping controller via particle swarm optimization (PSO) and an artificial bee colony (ABC) algorithm for doubly fed induction generator (DFIG) wind turbines. They aimed at MPPT and to decrease transient loads by controlling power transferred between the generator and the electrical grid in the presence of uncertainty. Aissaoui *et al.*, [8], presented in their work a comprehensive model of the WECS equipped with PMSG; they designed a Fuzzy-PI controller to maximize the extracted power with low power fluctuation. They managed to control and regulate the generator speed with low fluctuations.

The use of nonlinear control systems is considered one of the prevalent control systems of wind energy conversion systems. Thomsen, [9], described and analyzed different nonlinear control techniques for power and rotor speed regulation of the wind turbine. Additionally, other nonlinear control methods were implemented including gain scheduling technique, feedback linearization, [10], and sliding mode control, [11].

The work by Shao *et al.*, [12] deals with the restitution of the wind turbine pitch actuator system by addressing the PI- and PID-based pitch control methods. They sought to enhance the control system by mitigating the effect of pitch delay perturbations on the wind turbine output power.

Robust control theory tackled the control problem of WECS due to its ability to deal with external disturbances and model uncertainties. Bakka and Karimi, [13], implemented a mixed $H_2$-$H_\infty$ control design for the WECS, based on state-feedback control. They successfully regulated the rotor rotational speed subject to disturbances in the gearbox and wind turbine tower. Muhando *et al.*, [14], described a design of multi-objective $H_\infty$ control of the WECS that incorporates a doubly-fed induction generator. They designed a controller to accomplish the dual purpose of energy capture optimization and alleviating the cyclic load against wind speed fluctuations.

Linear varying parameter (LPV) controller is convenient for control problems that involve the regulation of wind turbine output power. Initially, the nonlinear model is transformed into an LPV model that consists of a group of linear models by assuming free-stream velocity, turbine shaft angular velocity, and pitch angle as varying parameters. Therefore, the control structure is reduced to an LPV controller which is a convex combination of linear controllers. Ying *et al.*, [15], have presented in their work a designed $H_\infty$ loop shaping torque and LPV-based pitch controllers. They aimed to enhance the performance of the pitch actuator system through the region of transition around the nominal wind speed. Gebraad *et al.*, [16], presented an LPV controller for WECS in the partial load region. They designed a full model that aimed to control the rotor vibration in the partial load region, and they used a proportional controller to enhance the produced power from the wind turbine. Inthamoussou *et al.*, [17], proposed an LPV controller for regulating the power of WECS above the rated wind speed. The suggested controller was compared with gain-scheduling PI and $H_\infty$ controllers. The work done by Lescher *et al.*, [18], adopted multivariable gain-scheduling controllers for the wind turbine based on a linear parameter-varying control approach. They designed a two-bladed wind turbine by situating smart micro-sensors hosted on the blades. They aimed to alleviate wind turbine cyclic loads addressed by the wind turbine in a full load regime.

The work presented here is mainly targeting the development of a robust linear parameter-varying $H_\infty$ controller for a WECS via pitch manipulation. It is precisely aimed at regulating generator output power via the regulation of the generator shaft angular velocity and torque. The control problem in hand lays down restrictions on the control design spec. Initially, the acceptable power fluctuations are limited to ± 5% of its nominal value regardless of the incoming wind speed variations. Additionally, the pitch actuator limitations introduce extra restrictions on the pitch angle and its derivative. These restrictions impose limitations on closed-loop performance. Furthermore, the suggested LPV controller is compared with the fuzzy logic controller under the same operating conditions.

## 2 Modelling of the Wind Turbine

The WECS operating region is dependent on the wind velocity, and it has three main regions: the No-Generation Region (NGR), the Partial-Load Region (PLR), and the Full-Load Region (FLR) as shown in Fig. 1.

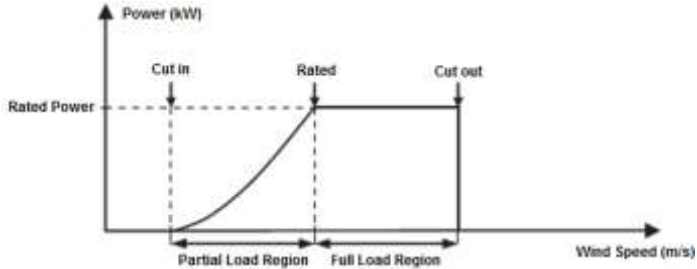

Fig. 1: Full operating region of variable speed WECS

The schematic diagram of the variable speed WECS is depicted in Fig. 2 (A), the figure shows four subsystems: the mechanical, the aerodynamic, the pitch actuator, and the generator subsystem. The wind-captured aerodynamic power ($P_r$) can be obtained by the following equation, [19];

$$P_r = \frac{1}{2}\rho\pi R^2 v^3 Cp(\lambda, \beta) \quad (1)$$

where $\rho$ is the air density, $R$ is the blade radius, and $v$ is the wind velocity. The aerodynamic torque ($T_r$) can be expressed by the following:

$$T_r = \frac{P_r}{\omega_r} \quad (2)$$

where $\omega_r$ is the rotor rotational speed. The power coefficient ($Cp$) can be given by Taher et al. [2]:

$$Cp(\beta, \lambda) = 0.22\left(\frac{116}{\lambda_i} - 0.6\beta - 5\right)\exp\left(-\frac{12.5}{\lambda_i}\right) \quad (3)$$

where:

$$\lambda_i = \cfrac{1}{\cfrac{1}{\lambda+0.12\beta} - \cfrac{0.035}{(1.5\beta)^3+1}} \quad (4)$$

where $\lambda$ is the tip speed ratio and $\beta$ is the pitch angle.

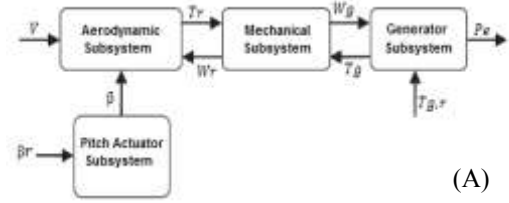

(A)

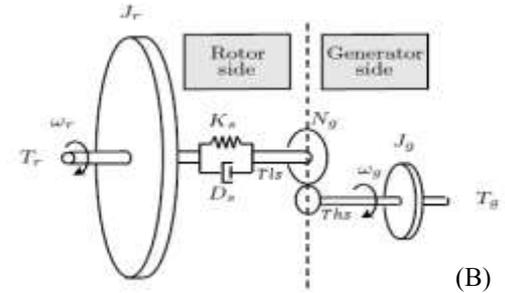

(B)

Fig. 2: Schematic diagram of the wind turbine (A), Two-mass WECS scheme (B).

The WECS model is considered a two-mass model as shown in Fig. 2(B). In this model, the turbine consists of two main components separated by the transmission: the low-speed shaft (rotor side) and the high-speed shaft (generator side). The gearbox ratio $N_g$ of the system is defined by:

$$N_g = \frac{T_{ls}}{T_{hs}} = \frac{w_g}{w_r} \quad (5)$$

where $T_{ls}$ is also defined by:

$$T_{ls} = D_s\dot{\delta} + K_s\delta \quad (6)$$

$$\dot{\delta} = \omega_r - \frac{\omega_g}{N_g} \quad (7)$$

where $\omega_g, T_{ls}, T_{hs}, \delta, K_s, D_s$, and $N_g$ are the generator speed, low-speed torque, high-speed torque, shaft twist, spring constant, damping coefficient, and gearbox ratio, respectively. The wind turbine mechanical dynamic equations [19] are obtained using Newton's second law:

$$\dot{\omega}_r J_r = T_r - T_{ls} \quad (8)$$

$$\dot{\omega}_g J_g = T_{hs} - T_g \quad (9)$$

where $J_r$ and $J_g$ are rotor inertia and the generator inertia. The generator model is given by:

$$\dot{T}_g = -\frac{1}{\tau_T}T_g + \frac{1}{\tau_T}T_{g,r} \quad (10)$$

where $T_{g,r}$ is the desired generator torque. This is a simplified generator first-order model with a time constant ($\tau_T$). The generator power can be acquired by:

$$P_e = T_g\,\omega_g \quad (11)$$

The pitch actuator is designed to regulate the rotational rotor speed at its rated value. It works by controlling the input power aerodynamic flow at the full load region. The turbine blades will turn in when the power is too low and will turn out when the power is too high. Generally, the power coefficient in (3) is minimized by raising the blades' pitch angle (β). The blade pitching process of the WECS imposes a time delay to reach the desired set-point value. Thus, the pitch actuator model is first-order with a rate limiter constrained to extreme values of ± 12 deg/s. The pitch actuator model shown in Fig. 3 is described as follows:

$$\dot{\beta} = \frac{1}{\tau_\beta}(\beta_r - \beta) \quad (12)$$

where $\beta_r$ and $\tau_\beta$ are the input blades' pitch angle and the time constant of the pitch actuator, respectively.

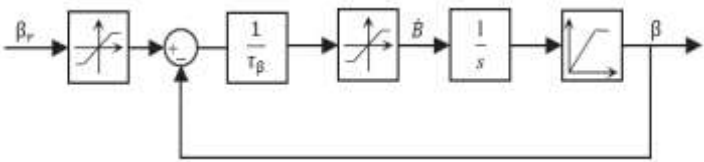

Fig. 3: Pitch actuator model block diagram

### 2.1 Nonlinear WECS
The equations of the mechanical dynamics obtained in Eq. (8) and (9) can be reformulated as:

$$\dot{\omega}_r = -\frac{K_s}{J_r}\delta - \frac{D_s}{J_r}\omega_r + \frac{D_s}{J_r N_g}\omega_g + \frac{1}{J_r}T_r \quad (13)$$

$$\dot{\omega}_g = \frac{K_s}{J_g N_g}\delta + \frac{D_s}{J_g N_g}\omega_r - \frac{D_s}{J_g N_g^2}\omega_g - \frac{1}{J_g}T_g \quad (14)$$

The state-space model for the non-linear wind turbine system is:

$$\dot{x} = Ax + Bu \quad (15)$$

$$\begin{bmatrix}\dot{\delta}\\\dot{w}_r\\\dot{w}_g\\\dot{\beta}\\\dot{T}_g\end{bmatrix} = \begin{bmatrix} w_r - \dfrac{w_g}{N_g} \\ \dfrac{D_s w_g}{J_r N_g} - \dfrac{D_s w_r}{J_r} - \dfrac{K_s \delta}{J_r} - \dfrac{T_r}{J_r} \\ \dfrac{D_s w_r}{J_g N_g} - \dfrac{D_s w_g}{J_g N_g^2} + \dfrac{K_s \delta}{J_g N_g} - \dfrac{T_g}{J_g} \\ -\dfrac{1}{\tau_\beta}\beta \\ -\dfrac{1}{\tau_T}T_g \end{bmatrix} + \begin{bmatrix}0 & 0\\0 & 0\\0 & 0\\ \dfrac{1}{\tau_\beta} & 0\\ 0 & \dfrac{1}{\tau_T}\end{bmatrix}\begin{bmatrix}\beta_r\\T_{g,r}\end{bmatrix} \quad (16)$$

where the state vector $x = [\delta\ w_r\ w_g\ \beta\ T_g]^T$, the control action $u = [\beta_r\ T_{g,r}]^T$ and the measured output $y = [w_g\ T_g]^T$.

### 2.2 Linearized WECS

The wind turbine aerodynamic torque ($T_r$) is a nonlinear function of rotor speed, wind speed, and pitch angle. The linearized model of the nonlinear system is realized via the first-order Taylor approximation approach. The nonlinear aerodynamic torque in Eq. (2), can be linearized as the following:

$$\hat{T}_r = K_w(\bar{\omega}_r, \bar{\beta}, \bar{v})\,\hat{\omega}_r + K_v(\bar{\omega}_r, \bar{\beta}, \bar{v})\,\hat{v} + K_\beta(\bar{\omega}_r, \bar{\beta}, \bar{v})\,\hat{\beta} \quad (17)$$

where $K_\omega$, $K_v$ and $K_\beta$ are the coefficients of linearization operating points defined as:

$$K_\omega = \left.\frac{\partial T_r}{\partial \omega_r}\right|_{(\bar{\omega}_r,\bar{\beta},\bar{v})},\quad K_v = \left.\frac{\partial T_r}{\partial v}\right|_{(\bar{\omega}_r,\bar{\beta},\bar{v})},\quad K_\beta = \left.\frac{\partial T_r}{\partial \beta}\right|_{(\bar{\omega}_r,\bar{\beta},\bar{v})}$$

.

The linearized rotor rotational speed can be obtained by substituting Eq. (17) in Eq. (13) to get the following formula:

$$\dot{\omega}_r = -\frac{K_s}{J_r}\delta + \frac{K_\omega - D_s}{J_r}\omega_r + \frac{D_s}{J_r N_g}\omega_g + \frac{K_v}{J_r}v + \frac{K_\beta}{J_r}\beta \quad (18)$$

The linear state-space model can be given by the following representation:

$$\Delta \dot{x} = A\Delta x + B\Delta u \quad (19)$$

$$\Delta y = C\Delta x \quad (20)$$

$$\begin{bmatrix}\Delta\dot{\delta}\\ \Delta\dot{\omega}_r\\ \Delta\dot{\omega}_g\\ \Delta\dot{\beta}\\ \Delta\dot{T}_g\end{bmatrix} = \begin{bmatrix} \omega_r - \frac{1}{N_g}\omega_g \\ -\frac{K_s}{J_r}\delta + \frac{K_\omega - D_s}{J_r}\omega_r + \frac{D_s}{J_r N_g}\omega_g + \frac{K_\beta}{J_r} \\ \frac{K_s}{J_g N_g}\delta + \frac{D_s}{J_g N_g}\omega_r - \frac{D_s}{J_g N_g^2}\omega_g - \frac{1}{J_g}T_g \\ -\frac{1}{\tau_\beta}\beta \\ -\frac{1}{\tau_T}T_g \end{bmatrix} + \begin{bmatrix}0\\ K_v/J_r\\ 0\\ 0\\ 0\end{bmatrix}v + \begin{bmatrix}0 & 0\\ 0 & 0\\ 0 & 0\\ \frac{1}{\tau_\beta} & 0\\ 0 & \frac{1}{\tau_T}\end{bmatrix}\begin{bmatrix}\Delta\beta_r\\ \Delta T_{g,r}\end{bmatrix} \quad (21)$$

where the state vector

$$\Delta x = x - \bar{x} = \{\Delta\delta\ \Delta\omega_r\ \Delta\omega_g\ \Delta\beta\ \Delta T_g\}^T,$$ the control action $\Delta u = u - \bar{u} = \{\Delta\beta_r\ \Delta T_{g,r}\}^T$, the measured output $\Delta y = y - \bar{y} = \{\Delta\omega_g\ \Delta T_g\}^T$ and $v$ is the exogenous input. The model transfer function is

$$G(s) = C(sI - A)^{-1}B \quad (22)$$

## 3 Controller Design

The robustness of a control system to disturbances has always been the main issue in feedback control systems. No need for feedback control if there are no disturbances in control systems, [20]. The robust control aims to achieve both robust performance and robust stability of the closed-loop system. In this section, LPV control based on mixed-sensitivity H∞ control is designed and presented above the rated wind energy conversion system. The primary objective of the proposed controller is to regulate the rotational rotor speed and the generator torque. Accordingly, maintaining the generator output power of the WECS to the rated power. Moreover, the turbulent wind velocity with large variations is presented in this research to estimate the stability and robustness of the suggested controller.

### 3.1 LPV Controller Design

The nonlinear WECS can be modeled as a linearized state-space system whose parameters vary with their states, [21], [22], [23]. The varying parameters ($\theta$) of the model are presumed to be measurable, bounded in a polytopic system, and slowly varying in real-time:

$$\dot{x} = A(\theta)x + B_1(\theta)v + B_2(\theta)u \quad (23)$$

$$y = C(\theta)x \quad (24)$$

The LPV model matrices are,

$$A(\theta) = \begin{bmatrix} \omega_r - \frac{1}{N_g}\omega_g \\ -\frac{K_s}{J_r}\delta + \frac{K_\omega(\theta) - D_s}{J_r}\omega_r + \frac{D_s}{J_r N_g}\omega_g + \frac{K_\beta(\theta)}{J_r} \\ \frac{K_s}{J_g N_g}\delta + \frac{D_s}{J_g N_g}\omega_r - \frac{D_s}{J_g N_g^2}\omega_g - \frac{1}{J_g}T_g \\ -\frac{1}{\tau_\beta}\beta \\ -\frac{1}{\tau_T}T_g \end{bmatrix},$$

$$B_1(\theta) = \begin{bmatrix}0 & \frac{K_v(\theta)}{J_r} & 0 & 0\end{bmatrix}^T,$$

$$B_2 = \begin{bmatrix}0 & 0 & 0 & \frac{1}{\tau_\beta} & 0\\ 0 & 0 & 0 & 0 & \frac{1}{\tau_T}\end{bmatrix}^T,$$

$$C = \begin{bmatrix}0 & 0 & 1 & 0 & 0\\ 0 & 0 & 0 & 0 & 1\end{bmatrix}.$$

The LPV controller is intended to control the WECS in the full-load regime, which covers wind speeds ranging from 11 to 24 m.s$^{-1}$. As a result, the wind speed (v) is the scheduling parameter. The rotor rotational speed $\omega_r$ of the wind turbine is held constant at a rated value of 4.3 rad/s. Where the main goal is to regulate the generator's output power around its rated value. The time-varying parameter varies in a polytopic system whose vertices are $\psi_1(v_{max})$ and $\psi_2(v_{min})$, in which the parameters are assumed to be measurable and slowly varying. The WECS LPV model can be realized with two vertices as follows:

$$\dot{x} = \left[\sum_{j=1}^{2} \alpha_j(\theta) A(\psi_j)\right] x + \left[\sum_{j=1}^{2} \alpha_j(\theta) B(\psi_j)\right] u \quad (25)$$

The controller $K(\theta)$ state space matrices at the vertices can be given by:

$$\begin{pmatrix} A_c(\theta) & B_c(\theta) \\ C_c(\theta) & D_c(\theta) \end{pmatrix} := \sum_{j=1}^{2} \alpha_j(\theta) K_j(\theta) = \sum_{j=1}^{2} \alpha_j(\theta) \begin{pmatrix} A_{cj}(\theta) & B_{cj}(\theta) \\ C_{cj}(\theta) & D_{cj}(\theta) \end{pmatrix} \quad (26)$$

where $\alpha_1 = \frac{v - v_{min}}{v_{max} - v_{min}}$, $\alpha_2 = 1 - \alpha_1$ and $\sum_{j=1}^{2} \alpha_j(\theta) = 1$.

According to Eq. (26), the LPV controller is designed as a convex combination of the vertices ψ1 and ψ2.

### 3.2 Mixed-weight H∞ Control Design

The mixed-weight H∞ controller's technique provides a closed-loop response of the system by shaping the frequency responses for noise attenuation and disturbance rejection. The controller design involves incorporating additional weighting functions in the original system, carefully chosen to demonstrate the system's performance and robustness specifications, [24]. As mentioned in Shaqarin et al., [25], [26], the generalized form of the mixed-weight H∞ problem can be elicited as explained in Fig. 4. Where w, y, z, u and e are the external input, the system-measured output, performance output, system input, and the control input, respectively. The general expanded plant $P(s)$ can be provided by:

$$z = \begin{bmatrix} z_1 \\ z_2 \\ z_3 \end{bmatrix} = \begin{bmatrix} W_1 e \\ W_2 u \\ W_3 y \end{bmatrix} \quad (27)$$

$$\begin{bmatrix} z \\ e \end{bmatrix} = P \begin{bmatrix} w \\ u \end{bmatrix} \rightarrow \begin{bmatrix} z \\ e \end{bmatrix} = \begin{bmatrix} P_{11} & P_{12} \\ P_{21} & P_{22} \end{bmatrix} \begin{bmatrix} w \\ u \end{bmatrix} \quad (28)$$

$$P = \begin{bmatrix} P_{11} & P_{12} \\ P_{21} & P_{22} \end{bmatrix} = \begin{bmatrix} W_1 & -W_1 G \\ 0 & W_2 \\ 0 & W_3 G \\ I & -G \end{bmatrix} \quad (29)$$

$$z = \begin{bmatrix} z_1 \\ z_2 \\ z_3 \\ e \end{bmatrix} = \begin{bmatrix} W_1 & -W_1 G \\ 0 & W_2 \\ 0 & W_3 G \\ I & -G \end{bmatrix} \begin{bmatrix} w \\ u \end{bmatrix} \quad (30)$$

where $W_1$, $W_2$ and $W_3$ are weighting functions. The closed-loop transfer function from $w$ to $z$ can be formulated using Eq. (28) and (29) as follows:

$$T_{zw} = P_{11}(s) + P_{12}(s) K(s) [I - P_{22}(s) K(s)]^{-1} P_{21}(s) \quad (31)$$

$$T_{zw} = \begin{bmatrix} W_1 S \\ W_2 KS \\ W_3 T \end{bmatrix} \quad (32)$$

where $S$ is the sensitivity function, $T$ is the complementary sensitivity function. The primary goal is to define a controller $(K)$ that reduces the infinity norm of $T_{zw}$ in the polytope $(\psi_1, \psi_2)$, such that $\|T_{zw}\|_\infty < \gamma$, where γ is the upper bound of $\|T_{zw}\|_\infty$. The selection of frequency-dependent weights ($W_1$, $W_2$, and $W_3$) substantially enhances the control design. Generally, at low frequencies, the sensitivity function $S$ is made small. This results in excellent disturbance rejection and a low tracking error. The complementary sensitivity function $T$ is also reduced in the high-frequency domain. As a result, there is high noise rejection and a broad stability margin. Further details on the weight selection can be found in Shaqarin et al. [25, 26].

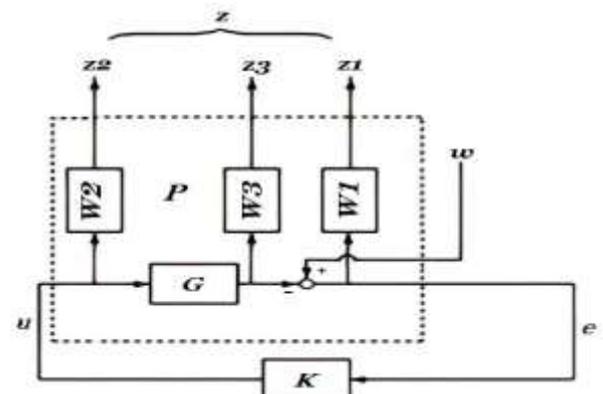

Fig. 4: Mixed-weight closed-loop system

### 3.3 PI-Fuzzy Logic Control (PIFLC) Design

The Fuzzy Logic Control (FLC) system has become one of the most common intelligent techniques utilized in many applications in current control systems. FLC can cope with various types of systems, ranging from linear processes to highly complex systems, such as nonlinear processes or time-varying systems. FLC versatility is attributed to its parameter tunability, such as the membership function type and number, rule base, scaling factor, inference techniques, fuzzification, and defuzzification. The process of fuzzy logic is shown in Fig. 5(A), which consists of a fuzzifier, inference engine, rule base, and defuzzifier. The process can be explained along these lines: the crisp inputs from the input data are initially fuzzified as fuzzy inputs. Hence, they trigger the inference engine and the rule base to generate the fuzzy output. The inference engine provides an input/output map just after blending the activated rules. Then, the inference engine outputs are sent to the defuzzifier, which produces the crisp outputs.

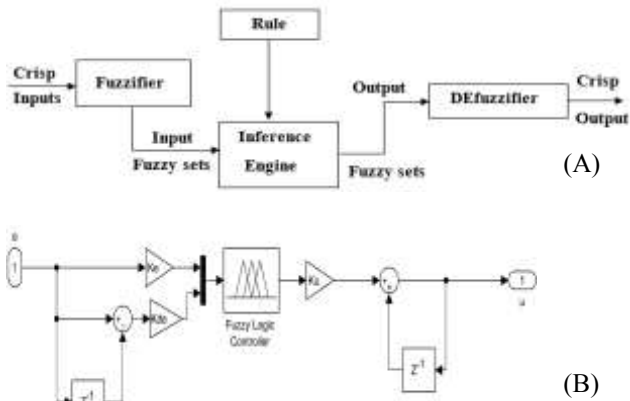

(A)

(B)

Fig. 5: Internal structure of the fuzzy controller(A), PIFLC control structure (B).

In the developed controller shown in Fig. 5(B), two input fuzzy variables: error (e) and change in error ($\Delta e$) with the output $\Delta u$ are shown. With a sampling period Ts, the signal $e$ is sampled and $\Delta e$ is calculated as:

$$\Delta e(k) = e(k) - e(k-1) \quad (33)$$

where $k$ denotes the sample number, and $z^{-1}$ denotes the unit time delay. As depicted in Fig. 5(B), the PIFLC controller output $u(k)$ can be found as:

$$u(k) = \Delta u(k) + u(k-1) \quad (34)$$

It is worth mentioning that the continuous control output $u(t)$ is obtained by assuming a zero-order hold between samples.

The membership functions of the inputs and the output are depicted in Fig. 6(A), where μ is the membership value. The performance of the PIFLC can be tuned via error gain ($K_e$), the change in the error ($K_{de}$), and the change in control output ($K_u$), as shown in Fig. 5(B).

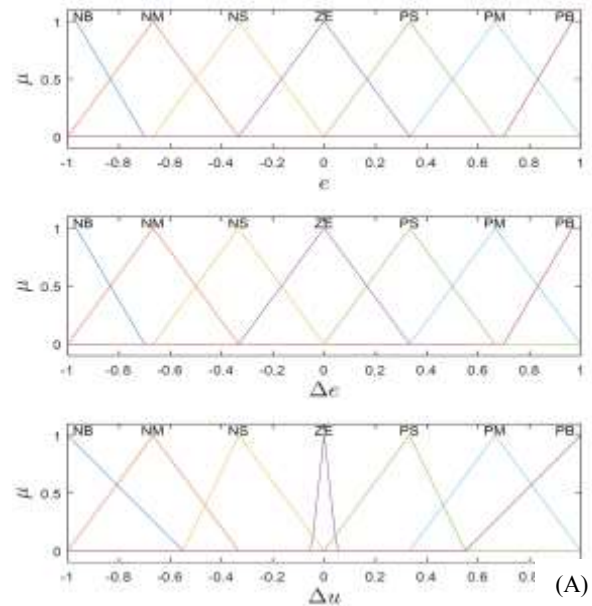

(A)

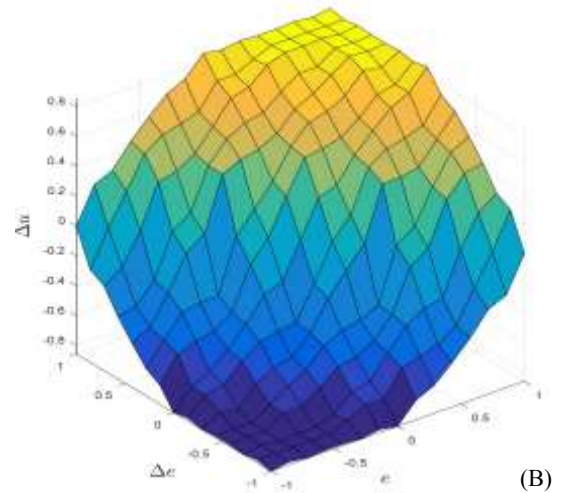

(B)

Fig. 6: Membership functions of the inputs e and $\Delta e$ and the output $\Delta u$ (A), Input-output surface relationship (B).

The FLC used in this work is a Mamdani type, where the structure of fuzzy rule is formulated as:

IF $e_f(k)$ is A AND $Ce_f(k)$ is B THEN $\Delta U(k)$ is C (33)

where $A, B$ and $C$ are fuzzy sets. The system has 49 fuzzy rules, and the surface input-output relationship is shown in Fig. 6 (B).

## 4 Results and Discussions

This section aims to assess the performance and stability of the suggested LPV-based $H_\infty$ control system of the WECS. This is accomplished via simulating the closed-loop response of the suggested controller using MATLAB/SIMULINK. Moreover, the suggested controller is compared with the response of the PIFLC and the WECS. In this work, the nominal parameters of the Vestas V29-225 kW WECS shown in Table 1 were used in the simulation. The linearization coefficient values $(K_\omega, K_v$ and $K_\beta)$ for both vertices are presented in Table 2. The weights are selected as follows: $W_1$ was chosen to yield better disturbance rejection through the shaping of the sensitivity function, which leads to a small tracking error. The weighting function $W_2$ was designed to shape the control sensitivity function, aiming at limiting the actuator effort. More precisely, $W_2^{-1}$ is responsible for limiting the pitch actuator effort to cope with the limited actuator bandwidth.

$$W_1 = \frac{s+10}{2s+0.1},$$

$$W_2 = \begin{bmatrix} \frac{s+60}{2000s+1.2e7} & 0 & 0 \\ 0 & \frac{s+2.5}{0.001s+25} & 0 \\ 0 & 0 & \frac{s+60}{2000s+1.2e7} \end{bmatrix}$$

$$W_3 = I$$

Table 1. Parameters of the Vestas WECS

| Parameters | Value |
|---|---|
| $P_{e0}$ | 225 kW |
| $\omega_{r0}$ | 4.3 rad/s |
| $\omega_{g0}$ | 105.78 rad/s |
| $\delta_0$ | 0.00655 rad |
| $J_g$ | 10 kg.m² |
| $J_r$ | 90000 kg.m² |
| $N_g$ | 24.6 |
| $R$ | 14.3 m |
| $D_s$ | 80000 s⁻¹ |
| $K_s$ | 8000000 N/m |
| $\tau_\beta$ | 0.15 s |
| $\tau_T$ | 0.1 s |
| $\rho$ | 1.225 kg.m3 |

Table 2. Coefficients of the linearization for the two vertices

| Vertex | $K_w(\bar{w}_r, \bar{\beta}, \bar{v})$ | $K_v(\bar{w}_r, \bar{\beta}, \bar{v})$ | $K_\beta(\bar{w}_r, \bar{\beta}, \bar{v})$ |
|---|---|---|---|
| $(v_{max})$ | $-6.57 \times 10^4$ | $1.618 \times 10^4$ | $-2.786 \times 10^4$ |
| $(v_{min})$ | $-6.91 \times 10^3$ | $1.229 \times 10^4$ | $-2.251 \times 10^3$ |

Using the above-mentioned control technique, the proposed LPV system based on an $H_\infty$ controller is intended to regulate the outputs of the WECS to the rated values without imposing large variations. The simulation of the closed-loop system of the WECS with both LPV controller and Fuzzy logic controller is discussed in the following sections.

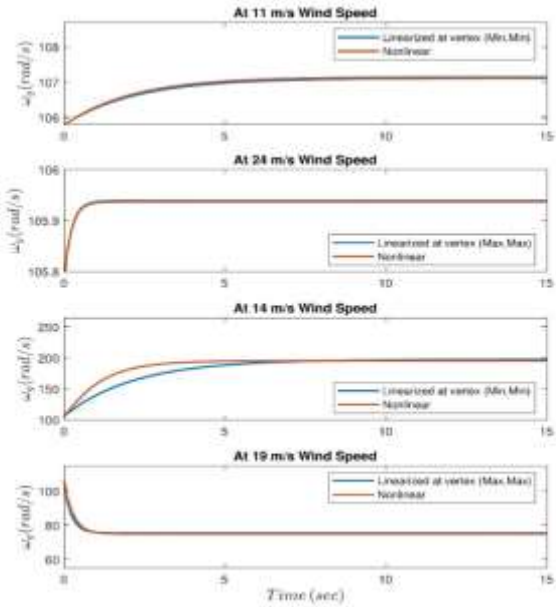

Fig. 7: Comparison of the nonlinear and linearized WECS open-loop step responses.

### 4.1 Nonlinear Versus Linearized WECS Simulation

The dynamic model of the WECS was obtained for both nonlinear and linearized cases as shown in sections (2.2 and 2.3). The linearization is carried out around operating points, which vary in a polytopic system whose vertices are

$\psi_1 (v_{max} = 24\ m/s,\ \beta_{max} = 24°)$ and

$\psi_2 (v_{min} = 11\ m/s,\ \beta_{min} = 0°)$. The rotational speed of the rotor ($\omega_r$) is assumed constant at the rated value.

Figure 7 shows the step response of the non-linear and the linearized systems under four wind speed values. The figure shows that the steady-state responses of the non-linear and linearized systems are identical at the linearization operating points with slight differences in the transients. However, the discrepancy between the two systems increases as the operating points change and move away from the linearization points.

### 4.2 Open-loop Response of WECS Subjected Wind Speed change

The open-loop response of WECS is simulated to evaluate the wind turbine performance when exposed to various wind speeds at a fixed pitch angle. The variable-speed wind turbine in Fig. 8 started with a smooth wind speed ranging from 11 m/s to 24 m/s at a minimum pitch angle of zero degrees. The figure shows that the rotor rotational speed and the speed of the generator increase as the free-stream velocity increases.

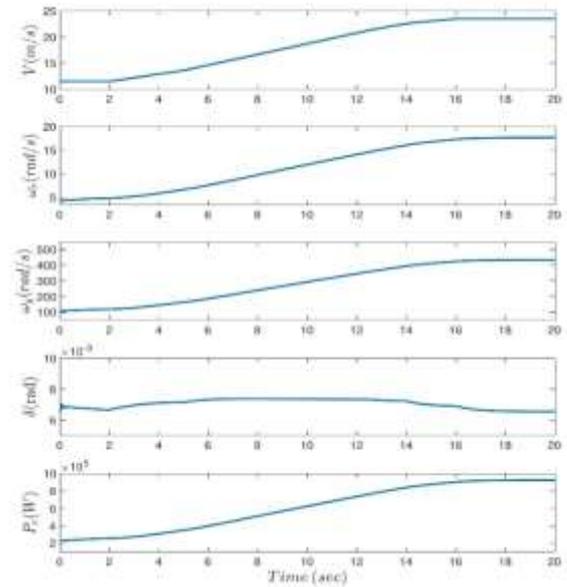

Fig. 8: Open-loop response of WECS subjected to a smoothly varying wind speed.

As a result, the generator's output power is increased to up to four times its nominal value. This motivates the need for closed-loop control of the WECS for power regulation. This is due to the fact that the generator's speed is also increased four times above its rated value, which Jeopardizes the wind turbine's safety and complicates the connection with the grid.

### 4.3 Controlled WECS Subjected to a Step-Ramp Change in Wind Speed

To evaluate the closed-loop system of the suggested LPV controller with variable speed WECS, step-ramp changes in free-stream velocity with a white noise variance of 0.0102 are introduced. These steps force the controller to modify the blade's pitch angle, and consequently the generator's rotational speed. Figure 9 illustrates the simulation of the WECS that started with a noisy free-stream velocity of 24 m/s, ranging from 0 to 25 sec. Then, the free-stream velocity stepped down with a slope of -0.6 m/s$^2$, until it reached a mean free-stream velocity of 17.5 m/s after 35 sec. This free-stream velocity was maintained constantly from t = 35 sec to 60 sec. Another stepped-down free-stream velocity occurred with the same previous slope until it reached the minimum free-stream of 11 m/s after 70 sec, then it remained constant. The figure shows that the response of the closed-loop system does not introduce high oscillations over the entire operating region.

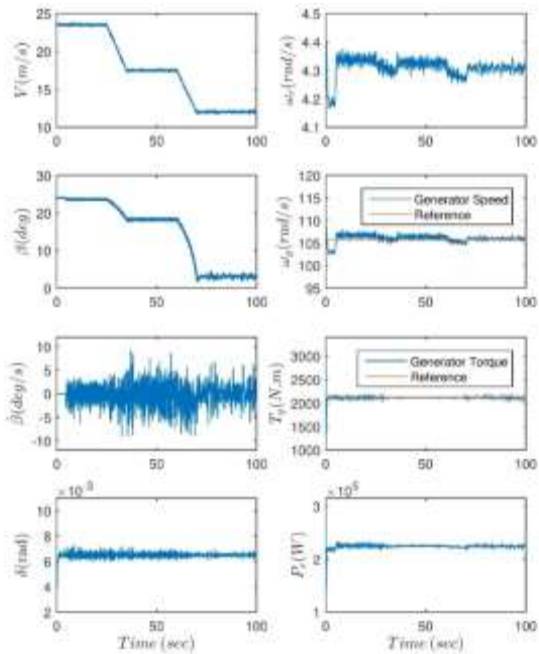

Fig. 9: Closed-loop simulation of the WECS with LPV controller subjected to step-ramp change in wind speed with white noise.

The response of the generator speed slightly changes around the nominal value, which is 105.78 rad/s. This indicates the LPV-based H∞ controller's effectiveness in maintaining the generator's output power very close to the rated value of 225 kW, in presence of noisy varying free-stream velocity, as seen in Fig. 9. Regarding the mechanical safety aspects, it is depicted in the figure that slight oscillations occurred in rotor shaft ($\pm\ 0.0005\ rad$), as seen in the response of the shaft twist angle over the whole operating range.

The suggested controller adequately maintained the generator speed and output power at their rated levels while maintaining the shaft's angle of twist nearly constant. It's worth noting that for the mechanical loads to be in the acceptable range, the less the peak twist angle, the less mechanical stress.

### 4.4 Closed-loop Response of the WECS with LPV and PI-Fuzzy Controllers in Presence of Turbulent Wind Speed

In this work, the Von Karman turbulence spectrum model is used with an average velocity of 17.5 m/s, with a turbulence intensity of the incoming wind flow greater than ± 10%, a turbulence length scale of 170 m, and a 2 m/s standard deviation. The wind speed used in the turbulence model shown in Fig. 10 varies in a range between 11 and 24 m/s, which lies in the full load region.

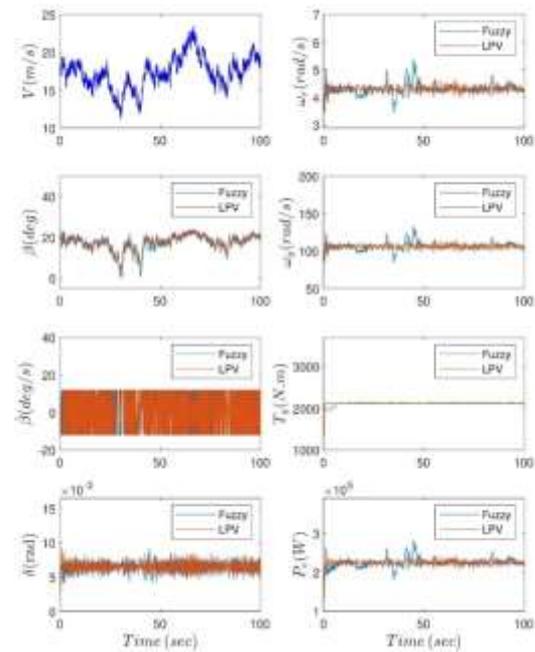

Fig. 10: Closed-loop simulation of the WECS with both LPV and PIFLC Controllers Subjected to Turbulent Wind Speed.

The proposed control system was able to handle these wind speed variations with high efficiency, as depicted in Fig. 10. The generator speed and the electromagnetic torque were controlled in the full load region of the WECS around their nominal values, which lead to stabilizing and maintaining the generator output power around its rated value, without violating the pitch angle ranges;

$0 < \beta < 24°$ and the pitch actuator constraints;

$-12 < \dot{\beta} < 12$. A PIFLC is designed and implemented in this paper to be compared with the proposed LPV control system. The gains were selected as presented in Table 3.

Table 3. PIFLC Gains of the pitch controller

| Gain | Value |
|---|---|
| $K_e$ | 2 |
| $K_{de}$ | 2 |
| $K_u$ | -2 |

The closed-loop simulations for both controllers are presented in Fig. 10. The figure shows the robustness of the proposed controller in enhancing the performance and stability of the WECS against free stream velocity and wind variations. The proposed controller response showed much fewer fluctuations in the generator outputs. This proves that the suggested LPV controller was superior at maintaining the output power of the generator at around its nominal value without causing large fluctuations in the output power of the WECS. On the other hand, there were significant power spikes and fluctuations in the Fuzzy controller response that exceeded the permissible design limitations. Though it was clear that the suggested LPV controller in this work was assigned to only one varying parameter (v) capable of satisfying the required control objectives.

The severe wind turbulence conditions with large mean wind speed variations implemented on the wind turbine, can undeniably conclude the main benefits of the LPV controller over the PIFLC controller, as shown in Fig. 10. The closed-loop response from the LPV controller indicates that the generator speed has very few fluctuations ($\pm 4\%$), whereas the PIFLC controller has significantly large peaks in the generator speed ($\pm 24\%$). This is translated to $\pm 4\%$ and $\pm 24\%$ peak fluctuations around the mean in the generated power for LPV and PIFLC controller cases, respectively. The fluctuations in the twist angle of the wind turbine shaft are comparable for both cases, whereas the variance of the fluctuation in the twist angle is two times less for the LPV controller case. This is highly beneficial in reducing the destructive mechanical loads on the wind turbine shaft. It is worth noting that the aforementioned analysis neglects the startup conditions.

## 5 Conclusion

In this paper, the suggested LPV based on an $H_\infty$ controller was employed to control a WECS via, manipulating the blades' pitch angle in the full load regime. The proposed controller was able to maintain and regulate the turbine shaft angular velocity, the electromagnetic torque, and thus the generator output power of the WECS to their nominal values. The proposed control design demonstrated proper performance and robustness when applied to a 225-kW WECS under turbulent free-stream velocity conditions. In comparison with the PIFLC, the suggested LPV controller was more effective in coping with the turbulent wind speed with a turbulence intensity of $\sim \pm 10\%$, which improved the wind turbine performance in terms of minimizing the fluctuations and smoothing the generator power. When the proposed LPV controller was assigned to only a single varying parameter (v), it showed its capability of meeting the desired control objectives of regulating and stabilizing the desired output power around its rated value, while complying with the pitch angle range; $0 < \beta < 24°$, and the pitch actuator constraints; $-12 < \dot{\beta} < 12$.